\documentclass[12pt]{article}
\usepackage{amsmath,amssymb,amsbsy,amsfonts,amsthm,latexsym,amsopn,amstext,
amsxtra,euscript,amscd}

  \newtheorem{thm}{Theorem}
  
  \newtheorem{lem}[thm]{Lemma}
  
  \theoremstyle{definition}
  
  \theoremstyle{remark}

\begin{document}

\def\squareforqed{\hbox{\rlap{$\sqcap$}$\sqcup$}}
\def\qed{\ifmmode\squareforqed\else{\unskip\nobreak\hfil
\penalty50\hskip1em\null\nobreak\hfil\squareforqed
\parfillskip=0pt\finalhyphendemerits=0\endgraf}\fi}

\def\cA{{\mathcal A}}
\def\cB{{\mathcal B}}
\def\cC{{\mathcal C}}
\def\cD{{\mathcal D}}
\def\cE{{\mathcal E}}
\def\cF{{\mathcal F}}
\def\cG{{\mathcal G}}
\def\cH{{\mathcal H}}
\def\cI{{\mathcal I}}
\def\cJ{{\mathcal J}}
\def\cK{{\mathcal K}}
\def\cL{{\mathcal L}}
\def\cM{{\mathcal M}}
\def\cN{{\mathcal N}}
\def\cO{{\mathcal O}}
\def\cP{{\mathcal P}}
\def\cQ{{\mathcal Q}}
\def\cR{{\mathcal R}}
\def\cS{{\mathcal S}}
\def\cT{{\mathcal T}}
\def\cU{{\mathcal U}}
\def\cV{{\mathcal V}}
\def\cW{{\mathcal W}}
\def\cX{{\mathcal X}}
\def\cY{{\mathcal Y}}
\def\cZ{{\mathcal Z}}

\def\MNL{{\mathfrak M}(N;K,L)}
\def\VNL{V_m(N;K,L)}
\def\RNL{R(N;K,L)}

\def\MNm{{\mathfrak M}_m(N;K)}
\def\VNm{V_m(N;K)}

\def\Xm{\cX_m}

\def \C {{\mathbb C}}
\def \F {{\mathbb F}}
\def \L {{\mathbb L}}
\def \K {{\mathbb K}}
\def \Q {{\mathbb Q}}
\def \Z {{\mathbb Z}}

\def\\{\cr}
\def\({\left(}
\def\){\right)}
\def\fl#1{\left\lfloor#1\right\rfloor}
\def\rf#1{\left\lceil#1\right\rceil}

\def \Prob{{\mathrm {}}}
\def\e{\mathbf{e}}
\def\em{\e_m}
\def\el{\e_\ell}
\def\Res{\mathrm{Res}}

\def\mand{\qquad \mbox{and} \qquad}

\newcommand{\comm}[1]{\marginpar{%
\vskip-\baselineskip 
\raggedright\footnotesize
\itshape\hrule\smallskip#1\par\smallskip\hrule}}

\title{On  Quadratic Fields Generated by Discriminants of
Irreducible  Trinomials}

\author{
{\sc Igor E.~Shparlinski\thanks{This work was supported in part by ARC Grant DP0556431}}\\ Department of Computing, 
Macquarie University\\
Sydney, NSW 2109, Australia\\
{\tt igor@ics.mq.edu.au}
}

\date{\today}

\maketitle

\begin{abstract}  A.~Mukhopadhyay, M.~R.~Murty
and K.~Srinivas have recently studied various arithmetic
properties of the discriminant  $\Delta_n(a,b)$ of the trinomial
$f_{n,a,b}(t) = t^n + at + b$, where $n \ge 5$ is  a fixed integer.
In particular, it is shown that, under the $abc$-conjecture,
for every $n \equiv 1 \pmod 4$,
the quadratic fields $\Q\(\sqrt{\Delta_n(a,b)}\)$ are pairwise distinct for
a positive proportion of such discriminants
with integers $a$ and $b$
such that $f_{n,a,b}$ is irreducible over $\Q$ and $|\Delta_n(a,b)|\le X$,
as $X\to \infty$.  We use the square-sieve and bounds of character
sums to obtain a weaker but unconditional version of this result.
\end{abstract}



\section{Introduction}

For a fixed integer $n\ge 2$, we use
$\Delta_n(a,b)$
to denote the discriminant of the trinomial
$$
f_{n,a,b}(t) = t^n + at + b.
$$

A.~Mukhopadhyay, M.~R.~Murty and K.~Srinivas~\cite{MMS}
have recently studied the arithmetic structure of $\Delta_n(a,b)$.
In particular, it is shown in~\cite{MMS}, under the $abc$-conjecture,
that if $n \equiv 1 \pmod 4$ then for a sufficiently large
positive $A$ and $B$ such that $B \ge A^{1+\delta}$ with some fixed
$\delta> 0$, there are at least $\gamma AB$ integers $a,b$ with
$$
A \le |a| \le 2A \qquad \text{and}\qquad
B \le |b| \le 2B
$$
and such that $f_{n,a,b}$ is irreducible and  $\Delta_n(a,b)$
is square-free, where $\gamma > 0$ depends only on $n$ and $\delta$.

Then this result is used  to derive
(still  under the $abc$-conjecture) that 
the quadratic fields
$\Q\(\sqrt{\Delta_n(a,b)}\)$ are pairwise distinct for
a positive proportion of such discriminants
with integers $a$ and $b$
such that $f_{n,a,b}$ is irreducible over $\Q$.

More precisely, for a real $X \ge 1$,
let $Q_n(X)$ be the number
of distinct fields $\Q\(\sqrt{\Delta_n(a,b)}\)$
taken for all of pairs of integers $a, b$ such that
$f_{n,a,b}$ is irreducible over $\Q$ and
$|\Delta_n(a,b)|\le X$.

Throughout the paper, we use  
$U = O(V)$, $U \ll V$, and $V \gg U$ 
as  equivalents of the inequality $|U| \le c V$ with some constant
$c> 0$, which  may depend only on  $n$.

It is shown in~\cite{MMS} that
for a fixed $n \equiv 1 \pmod 4$,
\begin{equation}
\label{eq:Q MMS bound}
Q(X) \gg X^{\kappa_n}, 
\end{equation}
where
$$ 
\kappa_n=\frac{1}{n} + \frac{1}{n-1}
$$
and $c_0 > 0$ is a constant depending only on $n$.  

It is also noted in~\cite{MMS} that the
Galois groups of irreducible 
trinomials $f_{n,a,b}$ have some interesting properties,
see also~\cite{CMS,HeSa,PlVi}. 
We remark that, since
\begin{equation}
\label{eq:Dab}
\Delta_n(a,b) = (n-1)^{n-1}a^n + n^n b^{n-1}
\end{equation}
for $n \equiv 1 \pmod 4$, there are
$O\(X^{\frac{1}{n} + \frac{1}{n-1}}\)$ integers $a$ and $b$ with
$|\Delta_n(a,b)|\le X$ and thus indeed~\eqref{eq:Q MMS bound}
means that
$$
Q(X) \gg \#\{(a,b) \in \Z^2~:~ |\Delta_n(a,b)|\le X\}.
$$

We use the square-sieve and bounds of character
sums to obtain a weaker but unconditional version of this result.
We note that, without the irreducibility of $f_{n,a,b}$ condition,
the problem of estimating $Q(X)$ can be viewed as
a bivariate analogue of the question, considered in~\cite{LuSh},
on the number of
distinct quadratic fields of the form
$\Q\(\sqrt{F(n)}\)$ for $n=M+1,\ldots, M+N$,
for a nonconstant polynomial $F(T) \in \Z[T]$.
Accordingly we use similar ideas, however we also exploit the
specific shape of the polynomial $\Delta_n(a,b)$
given by~\eqref{eq:Dab}.

\section{Main result}

In fact as in~\cite{LuSh} we consider a more general quantity than
$Q(X)$. Namely, for real positive $A$, $B$, $C$ and $D$ and a square-free integer
$s$, we denote by
$T_n(A,B,C,D;s)$ the
number of  pairs of integers
$$
(a,b) \in [C,C+A]\times [D,D+B]
$$
such that   $\Delta_n(a,b) = sr^2$ for some integer $r$. 

We write $\log x$ for the maximum of the natural 
logarithm of $x$ and 1, thus we 
always have  $\log  x \ge 1$.

\begin{thm}
\label{thm: TABCDs}
For real  $A\ge 1$, $B\ge 1$, $C\ge 0$ and $D\ge 0$ and a square-free $s$, 
we have
\begin{eqnarray*}
T_n(A,B,C,D;s) 
&\ll &(AB)^{2/3} \log (AB) +   A  \log (AB) + B  \log (AB)\\
 & & \qquad\qquad\qquad 
 +~(AB)^{1/3}  \(\frac{\log (ABCD) \log (AB)}{\log \log (ABCD)}\)^2. 
\end{eqnarray*}.
\end{thm}

Now, for real $A$, $B$, $C$ and $D$ we denote by
$S_n(A,B,C,D)$ the
number of distinct quadratic fields $\Q\(\sqrt{\Delta_n(a,b)}\)$
taken for all pairs of integers
$$
(a,b) \in [C,C+A]\times [D,D+B]
$$
such that $f_{n,a,b}$ is irreducible over $\Q$.
Using that the bound of Theorem~\ref{thm: TABCDs} is 
uniform in $s$, we derive

\begin{thm}
\label{thm: SABCD}
For real $A\ge 1$, $B\ge 1$, $C\ge 0$ and $D\ge 0$ ,
we have
\begin{equation*}
\begin{split}
S_n(A,B,C,D)  \gg
\min\Biggl\{\frac{(AB)^{1/3}}{\log (AB)},\,  \frac{A}{\log (AB)},\, &
 \frac{B}{\log (AB)} ,\\
 (AB)^{2/3} &  \(\frac{\log \log (ABCD)}{\log (ABCD) \log (AB)}\)^2\Biggr\}. 
\end{split}
\end{equation*}.
\end{thm}

The results of Theorems~\ref{thm: TABCDs} and~\ref{thm: SABCD}
are nontrivial in a very wide range of parameters $A$, $B$, $C$ and $D$
and apply to very short intervals.
In particular, $AB$ could be logarithmically small compared to $CD$.
Furthermore, taking
$$
A = C = \frac{1}{4 (n-1)^{n-1}}X^{1/n}
\qquad \text{and} \qquad
B = D = \frac{1}{4n^n}X^{1/(n-1)}
$$
  we see that
$$
Q(X)  \gg  X^{\kappa_n/3} (\log X)^{-1},
$$
which, although is weaker than~\eqref{eq:Q MMS bound},
does not depend on any unproven conjectures.

\section{Character Sums with the Discriminant}

Our proofs rest on some bounds for character sums. For an odd
integer $m$ we use $(w/m)$ to denote, as usual, the Jacobi symbol of
$w$ modulo $m$. We also put
$$
\em(w)= \exp(2 \pi i w/m).
$$

Given an odd integer $m\ge 3$ and arbitrary integers $\lambda, \mu$,
we consider the double character sums
$$
S_n(m; \lambda, \mu) =
\sum_{u,v=1}^{m}
\(\frac{\Delta_n(u,v)}{m}\) \em\(\lambda u + \mu v\).
$$

We need bounds of these sums in the case of $m =\ell_1\ell_2$
being a product of two primes $\ell_1> \ell_2\ge n$. 
However, using the multiplicative property of character sums
(see~\cite[Equation~(12.21)]{IwKow} for single sums, 
double sums behaves exactly the same way) we see that it is 
enough to estimate $S_n(\ell; \lambda, \mu)$ for primes $\ell$.

We start with evaluating these sums in the special case of $\lambda = \mu=0$
where we define
$$
S_n(\ell) = S_n(\ell; 0, 0) .
$$

\begin{lem}
\label{lem:Sm00} For $n \equiv 1 \pmod 4$ and a prime  $\ell$, we have
$$
S_n(\ell) \ll \ell.
$$
\end{lem}

\begin{proof} We can certainly assume that $\ell\ge n$ as otherwise the bound 
is trivial.

Recalling~\eqref{eq:Dab}, we derive
\begin{eqnarray*}
S_n(\ell) & = &
\sum_{u,v=1}^{\ell}
\(\frac{(n-1)^{n-1}u^n + n^n v^{n-1}}{\ell}\) \\
& = &  \sum_{u, v=1}^{\ell-1}  
\(\frac{(n-1)^{n-1}u^n + n^n v^{n-1}}{\ell}\) + O(\ell). 
\end{eqnarray*}
Substituting $uv$ instead of $u$, we obtain
\begin{eqnarray*}
S_n(\ell)
& = &  \sum_{v=1}^{\ell-1} \sum_{u=1}^{\ell-1}
\(\frac{(n-1)^{n-1}(uv)^n + n^n v^{n-1}}{\ell}\) + O(\ell) \\
& = &  \sum_{u, v=1}^{\ell-1}
\(\frac{\( (n-1)^{n-1}u^nv + n^n\) v^{n-1}}{\ell}\) + O(\ell)\\
& = &   \sum_{u, v=1}^{\ell-1}
\(\frac{(n-1)^{n-1}u^nv + n^n }{\ell}\)+ O(\ell) 
\end{eqnarray*}
since $n-1$ is even.
We now rewrite it in a slightly more convenient form as
$$
S_n(\ell) = \sum_{u=1}^{\ell-1} \sum_{v=1}^{\ell}
\(\frac{(n-1)^{n-1}u^nv + n^n }{\ell}\)+ O(\ell).
$$
As  $\gcd(\ell,n-1)=1$, making the
change of variables $(n-1)^{n-1}u^nv +  n^n= w$, we  note that
for every $u= 1, \ldots, \ell-1$ if 
$v =1, \ldots, \ell$, then $w$ runs through 
the complete residue system modulo $\ell$.
Hence,    
$$
S_n(\ell)=
(\ell -1) \sum_{w=1}^{\ell}
\(\frac{w  }{\ell}\) + O(\ell) = O(\ell),
$$
which concludes the proof.
\end{proof}

The following result can be derived from~\cite[Theorem~1.1]{FK}, 
however we give a self-contained and more elementary proof.

\begin{lem}
\label{lem:Sm arb} For $n \equiv 1 \pmod 4$, a  prime $\ell$ and 
arbitrary integers $\lambda,\mu$ with $\gcd(\lambda,\mu, \ell)=1$, we have
$$
|S_n(\ell;\lambda, \mu)| \ll \ell.
$$
\end{lem}

\begin{proof} As in the proof of Lemma~\ref{lem:Sm00}, 
we can certainly assume that $\ell\ge n$ as otherwise the bound 
is trivial.

Also as in the proof of Lemma~\ref{lem:Sm00},  we obtain 
$$
S_n(\ell;\lambda, \mu) = \sum_{u=1}^{\ell-1} \sum_{v=1}^{\ell}
\(\frac{(n-1)^{n-1}u^nv + n^n }{\ell}\)\el\((\lambda u + \mu) v\) + O(\ell).
$$
As $\gcd(\ell,n-1)=1$, making the
change of variables $(n-1)^{n-1}u^nv +  n^n = w$, we obtain
\begin{eqnarray*}
S_n(\ell;\lambda, \mu) & = &\sum_{u=1}^{\ell-1} \sum_{w=1}^{\ell}
\(\frac{w}{\ell}\)\el\((n-1)^{-n+1}u^{-n}(\lambda u + \mu) (w -n^n)\) + O(\ell)\\
& = &\sum_{u=1}^{\ell-1} 
\el\(-(n-1)^{-n+1}n^nu^{-n}(\lambda u + \mu)  \) \\
& &\qquad \qquad \quad  \sum_{w=1}^{\ell}
\(\frac{w}{\ell}\)\el\((n-1)^{-n+1}u^{-n}(\lambda u + \mu) w\) + O(\ell).
\end{eqnarray*}

The sums over $w$ is the Gauss sum, thus
\begin{eqnarray*}
\lefteqn{
\sum_{w=1}^{\ell}
\(\frac{w}{\ell}\)\el\((n-1)^{-n+1}u^{-n}(\lambda u + \mu) w\) }\\
& & \qquad \qquad \qquad \qquad \qquad = \(\frac{(n-1)^{-n+1}u^{-n}(\lambda u + \mu)}{\ell}\) \vartheta_\ell \ell^{1/2}, 
\end{eqnarray*}
for some complex $ \vartheta_\ell$ with $| \vartheta_\ell|=1$ (which
depends only on the residue class of $\ell$ modulo $4$), 
we refer to~\cite{IwKow,LN} for details. 

Since $n\equiv 1 \pmod 4$ we have
$$
\(\frac{u^{-n}}{\ell}\) = \(\frac{u^{-1}}{\ell}\) .
$$
Thus, combining the above identities we obtain
\begin{eqnarray*}
S_n(\ell;\lambda, \mu) & = & \vartheta_\ell \ell^{1/2}\sum_{u=1}^{\ell-1} 
 \(\frac{(n-1)^{-n+1}(\lambda  + \mu u^{-1})}{\ell}\)\\
& &\qquad \qquad \qquad  \el\(-(n-1)^{-n+1}n^nu^{-n}(\lambda u + \mu)  \) + O(\ell).
\end{eqnarray*}
Since $\gcd(\lambda,\mu, \ell)=1$, the Weil bound (see~\cite[Bound~(12.23)]{IwKow})
applies and implies that the sum over $u$ is $O(\ell^{1/2})$ which 
concludes the proof. 
\end{proof}

Combining Lemmas~\ref{lem:Sm00} and~\ref{lem:Sm arb}, and using the 
aforementioned multiplicativity property, we obtain

\begin{lem}
\label{lem:Sm compl} For $n \equiv 1 \pmod 4$, an  integer $m = \ell_1\ell_2$
which is a product of two distinct primes $\ell_1 \ne  \ell_2$ and arbitrary 
integers $\lambda,\mu$, we have
$$
|S_n(m; \lambda, \mu)| \ll m.
$$
\end{lem}

Finally, using the standard reduction
between complete and incomplete sums (see~\cite[Section~12.2]{IwKow}), 
we derive from Lemma~\ref{lem:Sm compl}

\begin{lem}
\label{lem:Sm incompl} For $n \equiv 1 \pmod 4$, 
 an   integer $m = \ell_1\ell_2$
which is a product of two distinct primes $\ell_1 \ne  \ell_2$ and 
real positive $A$, $B$, $C$ and $D$, we have
$$
\sum_{C \le a \le C+A}\,
\sum_{D \le b \le D+B}
\(\frac{\Delta_n(a,b)}{m}\)  \ll \(\frac{A}{m} + 1\)\(\frac{B}{m}+1\) m \log m.
$$
\end{lem}

\section{Irreducibility}

As in~\cite{MMS} we recall a very special case of
a results of S.~D.~Cohen~\cite{Coh} about
the distribution of irreducible polynomials
over a finite field $\F_q$ of $q$ elements.

\begin{lem}
\label{lem:Irred Trinom} For any prime $p$,  there are
$p^2/n + O(p^{3/2})$ irreducible trinomials
$t^n +\alpha t + \beta \in \F_p[t]$.
\end{lem}

\section{Proof of Theorem~\ref{thm: TABCDs}}

 For a real number $z \ge 1$ we let $\cL_z$ be the set of primes
$\ell \in [z, 2z]$.  For a positive integer $k$ we write
$\omega(k)$ for the number of prime factors of
$k$.

 We note that if
$k\ge 1$ is a perfect square, then for $z \ge 3$,
$$
\sum_{\ell \in \cL_z} \(\frac{k}{\ell}\) \ge  \# \cL_z - \omega(k),
$$

For each pair $(a,b)$, counted in  $T_n(A,B,C,D;s)$, we see t
hat $s\Delta_n(a,b)$ is a perfect square
and that $s\mid \Delta_n(a,b)$. Hence,
$$
\omega(s\Delta_n(a,b)) = \omega(\Delta_n(a,b)).
$$
Thus, for such $(a,b)$ we have
$$
\sum_{\ell \in \cL_z} \(\frac{s\Delta_n(a,b)}{\ell}\) \ge  \# \cL_z -
\omega_z(s\Delta_n(a,b)) =  \# \cL_z - \omega(\Delta_n(a,b)) .
$$
Since $\omega(k)! \le k$, we see from the Stirling formula
that 
$$
\omega(k) \ll \frac{\log k}{\log \log k}.
$$ 
Thus 
$$
\omega(\Delta_n(a,b))  \ll  \frac{\log (A+B+C+D)}{\log \log (A+B+C+D)}
\ll  \frac{\log (ABCD)}{\log \log (ABCD)}.
$$
In particular, by the Cauchy inequality,
\begin{eqnarray*}
\lefteqn{\(\# \cL_z\)^2T_n(A,B,C,D;s) }\\
  & & \qquad \ll \sum_{C \le a \le C+A}\,
\sum_{D \le b \le D+B}
\(\sum_{\ell \in \cL_z} \(\frac{s\Delta_n(a,b)}{\ell}\) + \omega(\Delta_n(a,b)) \)^2\\
     & & \qquad \ll   \sum_{C \le a \le C+A}\,
\sum_{D \le b \le D+B}
\(\sum_{\ell \in \cL_z} \(\frac{s\Delta_n(a,b)}{\ell}\) \)^2  \\
 & & \qquad \qquad \qquad\qquad\qquad\qquad\qquad
 +~AB  \( \frac{\log (ABCD)}{\log \log (ABCD)} \)^2. 
\end{eqnarray*}
We note that
\begin{eqnarray*}
\lefteqn{
\sum_{C \le a \le C+A}\,
\sum_{D \le b \le D+B}
\(\sum_{\ell \in \cL_z} \(\frac{s\Delta_n(a,b)}{\ell}\) \)^2 }\\
& & \qquad  = \sum_{C \le a \le C+A}\,
\sum_{D \le b \le D+B}
\(\(\frac{s}{\ell}\) \sum_{\ell \in \cL_z} \(\frac{\Delta_n(a,b)}{\ell}\) \)^2\\
& & \qquad \le \sum_{C \le a \le C+A}\,
\sum_{D \le b \le D+B}
\(\sum_{\ell \in \cL_z} \(\frac{\Delta_n(a,b)}{\ell}\) \)^2 .
\end{eqnarray*}
Squaring out and changing the order of summation, 
we obtain 
\begin{eqnarray*}
\lefteqn{\(\# \cL_z\)^2T_n(A,B,C,D;s) }\\
     & & \qquad \ll \sum_{\ell_1, \ell_2 \in \cL_z}  \sum_{C \le a \le C+A}\,
\sum_{D \le b \le D+B} \(\frac{\Delta_n(a,b)}{\ell_1\ell_2}\)  \\
 & & \qquad \qquad \qquad\qquad\qquad+  AB  \( \frac{\log (ABCD)}{\log \log (ABCD)}\)^2. 
\end{eqnarray*}
We now estimate the double sum over $a$ and $b$ trivially  as $O(AB)$
on the ``diagonal'' $\ell_1 = \ell_2$ and  use Lemma~\ref{lem:Sm incompl} 
otherwise, getting
\begin{equation}
\label{eq:Prelim}
\begin{split}
\(\# \cL_z\)^2T_n(A,B,C,D;s) 
 \ll  \# \cL_z AB +  \(\# \cL_z\)^2 &
\(\frac{A}{z^2} + 1\)\(\frac{B}{z^2}+1\) z^2 \log z\\
 +  AB & \( \frac{\log (ABCD)}{\log \log (ABCD)}\)^2. 
\end{split}
\end{equation}

By the prime number theorem we have $\# \cL_z \gg z/\log z$
so we derive from~\eqref{eq:Prelim} that
\begin{eqnarray*}
T_n(A,B,C,D;s) 
& \ll & ABz^{-1} \log z + AB z^{-2} + A \log z + B \log z\\
 & & \qquad  \qquad  
 +~z^2 \log z+ABz^{-2}  \( \frac{\log (ABCD) \log z}{\log \log (ABCD)}\)^2. 
\end{eqnarray*}

Clear the first term always dominates the second one, so the second term can 
be simply dropped. Thus taking $z = (AB)^{1/3}$ to balance the terms  
$ABz^{-1} \log z $ and $ z^2 \log z $, we obtain the desired estimate.

\section{Proof of Theorem~\ref{thm: SABCD}}

Let $p_0$ be smallest prime for which there exists 
an irreducible trinomial 
$$
t^n +\alpha_0 t + \beta_0 \in \F_{p_0}[t]
$$
($p_0$ exists by  Lemma~\ref{lem:Irred Trinom}). 

We now define the sets of integers
\begin{equation}
\label{eq:Sets AB}
\begin{split}
\cA & =~\{a \in [C,C+A]\cap \Z~:~ a \equiv \alpha_0 \pmod {p_0}\};\\
\cB & =~\{b \in [D,D+B]\cap \Z~:~ b \equiv \beta_0 \pmod {p_0}\}
\end{split}
\end{equation}
Clearly
\begin{equation}
\label{eq:CardAB}
\begin{split}
\# \cA \gg A \mand \# \cB \gg B  
\end{split}
\end{equation}
and every trinomials $t^n + at + b$ with $a\in \cA$, $b \in \cB$ is 
irreducible over $\Z$.  

Using Theorem~\ref{thm: TABCDs} to estimate the number 
of pairs $(a,b) \in  \cA\times \cB$ for which 
$\Q\(\sqrt{\Delta_n(a,b)}\)$ is a given quadratic field, we
obtain the desired result.

\section{Remarks}

Similar results can be obtained for more general trinomials
$t^n +at^m + b$ with fixed integers $n > m \ge 1$. Some properties
of the Galois group of these trinomials have been studied in~\cite{CMS,HeSa,PlVi}
where one can also find an explicit formula for their discriminant
(which generalises~\eqref{eq:Dab}).  In the case of $a=b=1$ it becomes 
$(-1)^{n(n-1)/2} \(n^n - (-1)^{n}m^m (n-m)^{n-m}\)$. Studying arithmetic 
properties of this expression, for example, its square-free part, 
when $n$ and $m$ vary in the region $N \ge n > m \ge 1$ for a
sufficiently large $N$, is a very challenging question.



\end{document}